\documentclass{amsart}
\usepackage{latexsym}
\newtheorem{theorem}{Theorem}[section]
\newtheorem{lemma}[theorem]{Lemma}
\newtheorem{proposition}[theorem]{Proposition}
\newtheorem{corollary}[theorem]{Corollary}
\theoremstyle{definition}                               
\theoremstyle{remark} 
\newtheorem*{remark}{Remark}
\numberwithin{equation}{section}
\begin{document}
\title[Automorphisms of the Weyl algebra]
{Automorphisms and ideals of the Weyl algebra}
%
%
	\author{Yuri Berest}
\address{Department of Mathematics, Cornell University,
         Ithaca, NY 14583-4201, USA}
\email{berest@math.cornell.edu}
%
%
%
%
%
%
     \author{George Wilson}
     \address{Department of Mathematics, 
      Imperial College, London SW7 2BZ, UK}
%
\email{g.wilson@ic.ac.uk}
%
   \thanks{The second author was supported in part by NSF Grant
DMS-94-00097; he wishes to thank MSRI and the University of 
California at Berkeley for their generous hospitality during 
the academic year 1998/9}
\begin{abstract} 
Let $A_1$ be the (first) Weyl algebra, and let $G$ be its 
automorphism group. We study the natural action of $G$ on the 
space of isomorphism classes of right ideals of $A_1$ (equivalently, 
of finitely generated rank 1 torsion-free right $ A_1$-modules). We show that
this space breaks up into a countable number of orbits each of 
which is a finite dimensional algebraic variety. Our results are 
strikingly
similar to those for the commutative algebra of polynomials in 
two variables; however, we do not know of any general principle that would 
allow us to predict this in advance. As a key step in the proof, 
we obtain a new description of the bispectral involution 
of \cite{W1}. We also make some comments 
on the group $G$ from the viewpoint of Shafaravich's theory of 
infinite dimensional algebraic groups.
\end{abstract}
\maketitle
\section{Introduction and statement of results}

The main aim of this paper is to give a simple description of the 
space $\mathcal{R}$ of isomorphism classes of (non-zero) 
right ideals in the complex Weyl algebra $A_1$, and of the action on 
$\mathcal{R}$ of the automorphism group of $A_1$. We begin by 
reviewing the corresponding (relatively trivial) results for 
the commutative algebra $A_0 = \mathbb{C}[x,y]$\,.

We call two non-zero ideals in $A_0$ {\it isomorphic} if they are 
isomorphic as $ A_0$-modules; equivalently, ideals $I$ and $J$ 
are isomorphic if we have $pI = qJ$ for some $p,q \in A_0$\,. 
Using the fact that reflexive ideals in $A_0$ are principal, 
it is easy to show that each isomorphism class contains a unique 
ideal of finite codimension. The ideals of a fixed codimension 
$n$ form a well-studied space, the Hilbert scheme 
$\mbox{\rm Hilb}_n (\mathbb{A}^2)$ of $n$-points in the affine 
plane. Alternatively, we can think of this space as the 
space of isomorphism classes of 
$n$-dimensional cyclic representations of $A_0$\,: to 
an ideal $I$ we assign the natural representation of $A_0$ 
on the quotient space $V = A_0/I$\,. A representation of $A_0$ 
on $V$ is nothing but a pair of commuting endomorphisms of $V$ 
(corresponding to the action of the generators $x$ and $y$ of $A_0$);
choosing a basis for $V$, we find that $I$ is represented by 
a pair of commuting matrices $(X,Y)$\,, uniquely determined up 
to simultaneous conjugation. Now let $\sigma$ be an automorphism 
of $A_0$\,: it induces an isomorphism (of rings) $\hat \sigma : A_0/I 
\to A_0/{\sigma(I)}$\,. If we use $\hat \sigma$ to identify these 
two spaces, and choose a basis, 
we find that the pair of matrices corresponding to 
$\sigma(I)$ is just $(\sigma^{-1}(X) , \, \sigma^{-1}(Y))$ (the 
notation means that we take the polynomials 
$\sigma^{-1}(x) , \, \sigma^{-1}(y) \in A_0$ and substitute 
$(X,Y)$ for $(x,y)$). To summarize:
\begin{enumerate}
\item The space of isomorphism classes of ideals breaks up 
into the disjoint union of the Hilbert schemes 
$\mbox{\rm Hilb}_n (\mathbb{A}^2)$\,.
\item A point of $\mbox{\rm Hilb}_n (\mathbb{A}^2)$ 
can be identified with a conjugacy class of pairs $(X,Y)$ of 
commuting $n \times n$ matrices (possessing a cyclic vector).
\item The formula for the action of an automorphism $\sigma$ 
on the pair of matrices $(X,Y)$ is the same as the 
formula for the action of $\sigma^{-1}$ on the generators 
$(x,y)$ of $A_0$\,.
\end{enumerate}
Note that, although the proofs of these statements are almost 
trivial, they all depend crucially on the first step, reducing 
the problem to the case of ideals of finite codimension.

Now let us turn to the Weyl algebra $A_1$\,. We recall that
$A_1$ is the associative algebra (over $\mathbb C$) generated 
by two elements $x$ and $y$ subject to the single relation 
$ [x,y] - 1 = 0 $.
In the main body of the paper we shall use the well known 
realization of $A_1$ as the algebra of ordinary differential 
operators with polynomial coefficients; in this Introduction 
we think of $A_1$ more abstractly, as a noncommutative 
version of $A_0$\,.
A fundamental difference between $A_0$ and $A_1$ is that 
$A_1$ has no finite dimensional representations (of positive 
dimension). Indeed, if $X$ and $Y$ are two $n \times n$ 
matrices with $[X,Y] - \mbox{\rm I} = 0$\,, then taking the trace, we 
obtain $n = 0$\,. Similarly, $A_1$  has no (left or right) ideals $I$ of  
finite codimension (for if it did, the quotient 
$A_1/I$ would be a finite dimensional representation). Thus, 
{\it a priori}, there seems no reason to expect that there 
should be results for $A_1$ similar to the ones indicated 
above for $A_0$\,: nevertheless, we shall see that this is the 
case. To obtain such results, we have to relax somewhat 
the notion of a representation, as follows. 
For each $n \geq 0$\,, let ${\mathcal C}_{n}$ be the space of
equivalence classes (modulo simultaneous conjugation) of
pairs $(X,Y)$ of $n \times n$ matrices satisfying the condition
\begin{equation}
\label{rk1}
[X, Y] - \mbox{\rm I} \ \, \mbox{\rm has rank at most}\ 1
\end{equation}
(the mysterious ``at most'' takes care of the case $n = 0$\,).
We might think of ${\mathcal C}_{n}$ as the space of 
``approximate $n$-dimensional representations'' of $A_1$. 
Let
$$
\mathcal{C} = \bigsqcup_{n \geq 0}\, {\mathcal C}_{n}
$$
be the disjoint union of the spaces ${\mathcal C}_{n}$.
In \cite{W2}, one of the authors constructed a bijective map
\smallskip
\begin{equation}
\label{1.3}
\beta : \mathcal C \to
\mbox{\rm Gr}^{\mbox{\rm \scriptsize ad}}
\end{equation}
from $\mathcal{C}$ to the adelic
Grassmannian that parametrizes rational solutions of the KP
hierarchy. This Grassmannian arises
(thinly disguised) also in the work of Cannings and Holland
on ideals in $A_1$\,. As for $A_0$\,, two 
right ideals $I$ and $J$ are called {\it isomorphic} if they are
isomorphic as right $A_1$-modules; equivalently, if we have
$pI = qJ$ for some $p, q \in A_1$.  Let $\mathcal R$ denote the
set of isomorphism classes of non-zero right ideals.
In \cite{CH}, Cannings and Holland constructed (in effect) 
a bijective map
\begin{equation}
\alpha : \mbox{\rm Gr}^{\mbox{\rm \scriptsize ad}}
\to \mathcal{R} \ .
\end{equation}
Set $\omega =  \alpha \beta$\,. Combining the results of \cite{CH}
and \cite{W2}, we obtain at once
\begin{theorem}
\label{th1}
The map
$\omega : \mathcal{C} \to \mathcal{R}$
is bijective.
\end{theorem}

Now let $G$ be the group of automorphisms of $A_1$. The natural
action of $G$ on $A_1$ induces an action on $\mathcal{R}$
(for if $I$ is a right ideal of $A_1$ and
$\sigma \in G$\,, then the isomorphism class of 
the right ideal $\sigma(I)$ clearly depends
only on that of $I$). We transfer this action to
$\mbox{\rm Gr}^{\mbox{\rm \scriptsize ad}}$ via the bijection
$\alpha$, and then to $\mathcal{C}$ via the bijection $\beta$.
To describe the resulting action of $G$ on $\mathcal{C}$,
recall (see \cite{D}, \cite{M2}) that $G$ is generated by the
special automorphisms of the form
$$
{\Phi}_{p}(D) = e^{p(y)}D\,e^{-p(y)}\,, \quad
{\Psi}_{q}(D) = e^{q(x)}D\,e^{-q(x)}\,,
$$
where $p$ and $q$ run over all polynomials (say without constant
term).
The action of these automorphisms on the generators $x$
and $y$ of $A_1$ is given by
\begin{equation}
\label{phieq}
\Phi_p(x) = x - p'(y)\ ,\quad \Phi_p(y) = y \ ;
\end{equation}
\begin{equation}
\label{psieq}
\Psi_q(x) = x\ ,\quad \Psi_q(y) = y + q'(x) \ .
\end{equation}
We shall prove
\begin{theorem}
\label{th2}
The action of $ G $ on $\mathcal{C}$ preserves each of the subspaces
$\mathcal{C}_{n}$\,. The action of the generators of $G$ on
$\mathcal{C}_{n}$ is given by the formulae
\begin{equation}
\label{phiEq}
\Phi_p(X,\, Y) = (X + p'(Y),\, Y) \,\, 
\end{equation}
\begin{equation}
\label{psiEq}
\Psi_q(X,\, Y) = (X,\, Y - q'(X)) \ .
\end{equation}
\end{theorem}
%

%
\begin{theorem}
\label{th3}
The action of $G$ on each of the spaces $\mathcal{C}_n$ is
transitive; that is, the $\mathcal{C}_n$ are the orbits of the
action of $G$ on $\mathcal{C}$.
\end{theorem}

Theorem~\ref{th3} is perhaps particularly unexpected, because 
it is a simpler result than we have in the commutative case.
Indeed, the automorphism group of $\mathbb{C}[x,y]$\,, or, 
equivalently, of the affine plane $\mathbb{A}^2$, does not 
act transitively on $\mbox{\rm Hilb}_n (\mathbb{A}^2)$ 
for $n > 1$\,, because any automorphism preserves the 
multiplicities of a collection of $n$ points of $\mathbb{A}^2$. 
See the last section of the paper for a further comment on this 
matter.

Of course, it follows at once from Theorem \ref{th3} that the spaces
$\beta(\mathcal{C}_n)$ are the orbits of the
action of $G$ on $\mbox{\rm Gr}^{\mbox{\rm \scriptsize ad}}$;
and that the spaces $\omega(\mathcal{C}_n)$ are the orbits of the
action of $G$ on $\mathcal{R}$\,. Each of these facts has
considerable independent interest. In \cite{W2} the spaces
$\beta(\mathcal{C}_n)$ were characterized as the union of all
the $n$-dimensional open cells in
$\mbox{\rm Gr}^{\mbox{\rm \scriptsize ad}}$; but our present
description of them as the orbits of a group action seems more
straightforward. The problem of describing the orbits of $G$ in
$\mathcal{R}$ goes back at least to Stafford's paper \cite{St2}:
it is known to be equivalent to the problem
of classifying rings of differential operators on affine curves;
or again, to classifying the algebras Morita equivalent to $A_1$ 
(see, for example, \cite{BW}, \cite{CH}, \cite{K}, \cite{St2}, 
\cite{W3}). 
The key fact that the $G$-orbits in $\mathcal{R}$
are classified by a single non-negative integer $n$ was first
discovered in the 1994 thesis \cite{K} of Kouakou: however, Kouakou's
work was not published at that time, and remained unknown to
us until after we rediscovered his main result in \cite{BW}.

Here are a few words about our proof of Theorem \ref{th2}.
The difficulty is that the maps $\alpha$ and $\beta$ do not respect
the symmetry between $x$ and $y$ (or $X$ and $Y$), so that although
the formula (\ref{phiEq}) is easy to check, the similar
formula (\ref{psiEq}) is not: indeed, we shall see that the automorphisms
$ \Phi_p $ act on $\mbox{\rm Gr}^{\mbox{\rm \scriptsize ad}}$ as the group
of KP flows, whereas we do not know any such simple description of the action
of the $ \Phi_q $.
To pinpoint the problem, we introduce the
{\it formal Fourier transform}\, $\varphi$\,: it is the automorphism of
$A_1$ whose action on the generators is given by
\begin{equation}
\varphi(x) = -y \ ,\ \ \varphi(y) = x \ .
\end{equation}
We have ${\Psi}_{p} = \varphi \,{\Phi}_{p} {\varphi}^{-1}$, hence the
$\Phi_p$ and $\varphi$ generate $G$\,; thus our problem is reduced
to tracking how $\varphi$ acts on
$\mbox{\rm Gr}^{\mbox{\rm \scriptsize ad}}$, and thence on
$\mathcal{C}$. Although the action of $\varphi$
on $\mbox{\rm Gr}^{\mbox{\rm \scriptsize ad}}$ does not seem
to have any simple description, we can write $\varphi$
as the composition $\varphi = bc$ of two involutory
{\it anti}-automorphisms
of $A_1$: the action of $b$ and $c$ on the generators of $A_1$
is defined by
\begin{equation}
\label{ab}
b(x) = y\ ,\  b(y) = x \ ;\quad c(x) = -x \ ,\ c(y) = y \ .
\end{equation}
(If we think of $A_1$ as an algebra of differential operators, 
then $c$ is the formal adjoint). Using the fact that
the ideals in $A_1$ are reflexive $A_1$-modules, we can extend the action
of $G$ on $\mathcal{R}$ to an action of the larger group
$\widehat G$ of automorphisms and anti-automorphisms of $A_1$;
we then find that the actions of $b$ and $c$ on
$\mbox{\rm Gr}^{\mbox{\rm \scriptsize ad}}$ correspond to known
involutions. In particular, we shall see
\begin{theorem}
\label{th4}
The anti-automorphism $b$ acts on
$\mbox{\rm Gr}^{\mbox{\rm \scriptsize ad}}$\!
as the bispectral involution introduced in \cite{W1}.
\end{theorem}

Theorem \ref{th4} is (to us) the most difficult result
in this paper: in the following Sections 2--8 we prepare
the machinery for its proof, and also summarize the most 
essential definitions and facts that we need from the 
literature. The new material in these sections concerns 
some spaces of differential operators which we denote by 
$ {\mathcal D}(U,V) $. Here $ U $ and $ V $ are certain 
linear subspaces of the space $ {\mathbb C}(z) $ of 
rational functions in one variable, and $ {\mathcal D}(U,V) 
\subset {\mathbb C}(z)[\partial/\partial z] $ is the set of 
all differential operators with rational coefficients that
map $ U $ into $ V $ (mostly $ U $ and $ V $ will represent 
points of the Grassmannian $\mbox{\rm Gr}^{\mbox{\rm \scriptsize ad}}$).
These spaces $ {\mathcal D}(U,V) $ play a basic role in the paper 
of Cannings and Holland, but (to our knowledge) they have not
so far appeared in the theory of integrable systems. In that 
theory, on the other hand, an essential part of the 
machinery involves the {\it Baker functions} \,$ \psi_{U} $
of points $ U \in \mbox{\rm Gr}^{\mbox{\rm \scriptsize ad}} $:
these functions do not appear in the work of Cannings and Holland. 
We bridge this gap in Section~8, which gives a description of 
$ {\mathcal D}(U,V) $ in terms of the Baker functions of $ U $ and 
$ V $. Our main theorems all follow
easily from this (and other facts available in the literature):
the proofs are given in Sections~9~and~10. Finally, in the last 
section of the paper we offer a preliminary attempt to 
interpret our results from the perspective of algebraic groups.
  
%
%
%
%

In connection with Theorem~\ref{th1}, there are at least 
two natural questions which remain untouched in the present paper: 
one is to find a more direct definition of the map $ \omega $ 
(without passing 
through the Grassmannian, where some of the symmetry is difficult 
to see), another is to describe explicitly its inverse. We shall 
address both these problems in our subsequent paper \cite{BW2}, 
where the correspondence $ \omega $ will be constructed by
cohomological methods developed recently in the 
framework of noncommutative algebraic geometry.
The idea of such an approach goes back to the work
of L.~Le Bruyn \cite{L}. In \cite{BW2} we hope also to clarify 
the relationship between our stratification of 
$ {\mathcal R} $ and Le Bruyn's ``moduli spaces''.


With the exception of Theorem \ref{th2}, the main results of this 
article were announced in \cite{W3}.

\section{The group $\Gamma$ and the spaces $\mathcal{D}(U,V)$}

In this section we confront the following technical problem:
a point $ W $ of $ \mbox{\rm Gr}^{\mbox{\rm \scriptsize ad}} $
is a (certain kind of) linear space of rational functions,
but we need to consider the action on 
$ \mbox{\rm Gr}^{\mbox{\rm \scriptsize ad}} $ of the group
$ \Gamma $ of KP flows. Loosely speaking, this action is given
by multiplication by exponential functions. The most natural
way to make sense of that is to work temporarily with the 
completion of $ W $ in some larger function space. Because
we need to consider also the action of $ \Gamma $ on the spaces
of differential operators $ {\mathcal D}(U,V) $ it is 
inconvenient to use $\mbox{\rm L}^2$ completions (as in \cite{W1}, 
\cite{W2}); instead, we shall work with holomorphic 
functions, adapting to our circumstances a suggestion of
\cite{ACKP}.

To begin with, we work in slightly greater generality
than will be needed for the main part of the paper.
Let $\mathcal{S}$ be the set of all linear subspaces
$V \subset \mathbb{C}[z]$ such that we have 
\begin{equation}
\label{Sdef}
\rho(z) \mathbb{C}[z] \subset V \ \text{for some polynomial}\  
\rho\ ;
\end{equation}
and let $\Gamma$ be the group of all functions of the form
$\gamma(z) = e^{p(z)}$\,, where $p$ is a polynomial. 
To define the action of $ \Gamma $ on $\mathcal{S}$, we introduce
the space $\mathcal{H}$ of entire functions
on $\mathbb{C}$\,, with its usual topology (uniform convergence
on compact subsets), and denote by 
$\overline{\mathcal S}$ the set 
of all closed subspaces $\mathcal{V} \subset \mathcal{H}$ such 
that we have 
$$
\rho(z) \mathcal{H} \subset \mathcal{V} \ 
\text{for some polynomial}\  \rho\,.
$$
If $V \in \mathcal{S}$\,, let $\overline V$ be the closure of
$V$ in $\mathcal{H}$\,; then  clearly, 
$\overline V \in \overline{\mathcal{S}}$\,.
If $ \mathcal{V} \in \overline{\mathcal{S}}$\,, we set  
$\mathcal{V}^{\mbox{\rm \scriptsize alg}} = \mathcal{V} \cap \mathbb{C}[z]$\,; 
clearly, $\mathcal{V}^{\mbox{\rm \scriptsize alg}} \in \mathcal{S}$\,. 
If $\rho$ is a polynomial of degree $n$\,, then we have the 
decomposition 
$$
\mathcal{H} = \rho \mathcal{H} \oplus \{ 
\text{polynomials of degree less than}\ n \}
$$
of $\mathcal{H}$ into a direct sum of two closed subspaces. Comparing 
this with the similar decomposition of $\mathbb{C}[z]$\,, we 
easily find
\begin{lemma}
The maps $V \mapsto \overline V$ and $\mathcal{V} \mapsto 
\mathcal{V}^{\mbox{\rm \scriptsize alg}}$ define inverse bijections 
between $\mathcal{S}$ and $\overline{\mathcal{S}}$\,.
\end{lemma}

We use this Lemma to transfer the obvious action of $\Gamma$ on 
$\overline{\mathcal{S}}$ to $\mathcal{S}$\,; that is, if 
$V \in \mathcal{S}$ and $\gamma \in \Gamma$\,, we define
$$
\gamma V = (\gamma \overline V)^{\mbox{\rm \scriptsize alg}} \ .
$$

As in the Introduction, if \,$U,\, V \in \mathcal{S}$\,, 
we write
\begin{equation}
\label{duv}
\mathcal{D}(U,V) = \{ D \in \mathbb{C}(z)[\partial] :
D.U \subset V \}
\end{equation}
for the set of all differential operators with {\it rational}
coefficients that map $U$ into $V$\,. Here
$\partial \equiv \partial / \partial z$\,; we write $D.f$
for the action of $D$ on the function $f$ (to avoid confusion
with the composition $Df$ of $D$ with the multiplication
operator $f$). If $\rho$ is a polynomial, then we clearly 
have 
$$
\mathcal{D}(\rho U,\rho V) = \rho \mathcal{D}(U,V) \rho^{-1}\ .
$$
We need to know that this is still true if we replace $\rho$ 
by an element of $\Gamma$\,. 
\begin{lemma}
\label{D}
Let $U,\, V \in \mathcal{S}$\,, $D \in \mathbb{C}(z)[\partial]$\,.
Then
$$
D.U \subset V \Longleftrightarrow D.{\overline U} \subset \overline{V}\ .
$$
\end{lemma}
\begin{proof}
The implication ``$\Leftarrow$'' is trivial. Conversely, suppose
$D.U \subset V$\,, and choose a polynomial $\pi$ such that $\pi D$
has polynomial coefficients. Then $\pi D.U \subset \pi V$\,; since
$\pi D$ is a bounded operator on $\mathcal{H}$\,, it follows that
$\pi D.{\overline U} \subset \overline{\pi V} = \pi \overline V$\,.
Hence $D.{\overline U} \subset \overline{V}$\,, as required.
\end{proof}
\begin{corollary}
\label{hardlem}
Let $U,\, V \in \mathcal{S}\,, \,\gamma \in \Gamma$\,. Then
$$
\mathcal{D}(\gamma U,\gamma V) = \gamma \mathcal{D}(U,V) \gamma^{-1}\ .
$$
\end{corollary}
\begin{proof}
The corresponding result for the closures $\overline{U}, \,\overline{V}$ 
is trivial, so the Corollary follows at once from Lemma \ref{D}
and the fact that the operators in 
$ \gamma \mathcal{D}(U,V) \gamma^{-1} $ have rational coefficients.
\end{proof}
\section{The Cannings-Holland correspondence}

For the rest of this paper we confine our attention to the spaces
$V \in \mathcal{S}$ that are (in the terminology of
\cite{CH}) {\it primary decomposable}. The definition is
as follows. First, if $\lambda \in \mathbb{C}$\,, a space
$V \in \mathcal{S}$ is called $\lambda$-{\it primary} if
the polynomial $\rho$ in (\ref{Sdef}) can be chosen to be a 
power of $z - \lambda$\,, that is, if
$V$ contains the
ideal $(z - \lambda)^r \mathbb{C}[z]$ for some $r \geq 0$\,.
Then $V$ is called {\it primary decomposable} if 
it is a finite intersection of $\lambda$-primary
subspaces (for various $\lambda$). The action of $\Gamma$ on
$\mathcal{S}$ preserves the class of primary decomposable
subspaces (for it clearly preserves the class of $\lambda$-primary
subspaces for any fixed $\lambda$).
Primary decomposable subspaces play a fundamental role
in both of the papers \cite{CH} and \cite{W1}: indeed, this 
observation was the starting point for the present work.

To formulate the results of Cannings and Holland,
we think of the Weyl algebra $A_1$ as the ring
${\mathbb C}[z,\, \partial / \partial z]$
of ordinary differential operators with polynomial coefficients.
If $V \in \mathcal{S}$\,, we   
assign to it the right ideal
$$
\alpha(V) = \mathcal{D}({\mathbb C}[z], V)
$$
of $A_1$\,. If $I$ is a right ideal in $A_1$, we assign to it the
linear subspace $ e(I) $ spanned by all polynomials 
of the form $ D.f $ with $ D \in I $, $ f \in {\mathbb C}[z] $.
Equivalently, we have simply
$$
e(I) = \lbrace D.1 : D \in I \rbrace \subset {\mathbb C}[z] \ .
$$
The following theorem is proved in \cite{CH}.
\begin{theorem}
\label{ch}
The maps $\alpha$ and $e$ define inverse bijections
between the set of primary decomposable subspaces of
${\mathbb C}[z]$ and the set of right ideals of $A_1$ that
have non-zero intersection with ${\mathbb C}[z]$\,.
\end{theorem}
Now, it is easy to show that 
every non-zero ideal in $A_1$ is isomorphic (as right $A_1$-module)
to one that intersects
$\mathbb C[z]$ (see \cite{St2}, Lemma 4.2).
Further, two ideals $I$ and
$J$ that intersect $\mathbb C[z]$ are isomorphic if and
only if $\pi(z)I = \rho(z)J$ for some polynomials $\pi$ and $\rho$~.
The maps $\alpha$ and $e$ commute with left
multiplication by polynomials, so two primary decomposable
subspaces $U$ and $V$ correspond to isomorphic ideals
if and only if \,$\pi(z)U = \rho(z)V$\, for some $\pi$ and $\rho$ . 
Thus Theorem \ref{ch} implies the following
\begin{corollary}
\label{chcor}
The maps $\alpha$ and $e$ determine inverse bijections
between the set $\mathcal{R}$ of isomorphism classes 
of non-zero right ideals of 
$A_1$ and the set of classes of primary decomposable
subspaces of $\mathbb C[z]$ modulo the equivalence 
relation 
\begin{equation}
\label{equiv}
U \equiv V \Longleftrightarrow \pi(z)U = \rho(z)V 
\ \text{\rm for some polynomials}\ \pi \ \text{\rm and} \ \rho \ .
\end{equation}
\end{corollary}
\section{Duality and the action of $\widehat G$ on $\mathcal{R}$}

Let $Q$ denote the quotient (skew) field of $A_1$. We recall
that a right $A_1$-submodule $I$ of $Q$ is called a
{\it fractional right ideal} if there is some $a \in A_1$
such that $aI \subset A_1$ (and hence $aI$ is a right
ideal in the usual sense). Fractional left ideals are defined
similarly. If $I$ is a fractional right ideal,
its {\it dual} $I'$ is defined by
$$
I' = \{q \in Q : qI \subset A_1\} \ .
$$
Clearly, $I'$ is a (fractional) left ideal; the terminology
''dual'' is justified by the easy
\begin{lemma}
As left $A_1$-module, $I'$ is isomorphic to
$\mbox{\rm Hom}_{A_1}(I, A_1)$\,.
\end{lemma}

Similarly, if $J$ is a (fractional) left ideal, the dual
right ideal is defined by
$$
J' = \{q \in Q : Jq \subset A_1\} \ .
$$
A basic property of $A_1$ (see \cite{MR}) is that it is a
{\it hereditary} Noetherian ring, which means that its ideals
are projective, and hence reflexive, as $A_1$-modules.
That gives us the important
\begin{lemma}
\label{reflem}
For any right ideal $I$ and left ideal $J$, we have
$$
(I')' = I\quad,\quad (J')' = J \ .
$$
\end{lemma}

Using Lemma \ref{reflem}, we can extend the obvious action of 
$G = \mbox{\rm Aut}(A_1)$ on $\mathcal{R}$ to an action of 
the larger group
$\widehat G$ of all automorphisms and anti-automorphisms of
$A_1$. Recall from the Introduction that if $\sigma \in G$
is an automorphism of $A_1$, then we make $\sigma$ act on
(the class of) the right ideal $I$ by
\begin{equation}
\label{act1}
\sigma \star I = \sigma(I) \equiv \{\sigma(i) : i \in I\}\ .
\end{equation}
If now $\sigma$ is an {\it anti}-automorphism of $A_1$ and $I$
is a right ideal, then the set $\sigma(I)$ is a {\it left}
ideal, so its dual is again a right ideal; we define the
action of $\sigma$ on (the class of) $I$ by
\begin{equation}
\label{act2}
\sigma \star I = \sigma(I)' \ \, \text{if $\sigma$ is an anti-automorphism}\ .
\end{equation}
\begin{lemma}
The formulae (\ref{act1}) and (\ref{act2}) define an action
of $\widehat G$ on $\mathcal{R}$\,.
\end{lemma}
\begin{proof}
Everything is trivial, except possibly to check that if
$\sigma$ and $\tau$ are anti-automorphisms, and
$\theta = \sigma \tau$\,, then $\sigma \star
(\tau \star I) = \theta \star I$\,. For this, we need to
prove that
$$
\sigma[\tau(I)']' = \sigma[\tau(I)] \ .
$$
But because of Lemma \ref{reflem}, it is equivalent to
prove that
$$
\sigma[\tau(I)'] = \sigma[\tau(I)]' \ ,
$$
which is easy.
\end{proof}
\section{The adelic Grassmannian}
There does not seem to be any natural way of choosing 
representatives for the equivalence classes of primary 
decomposable subspaces in Corollary \ref{chcor}. However, 
these classes are in 1-1 correspondence with the points of 
the {\it adelic Grassmannian} 
$\mbox{\rm Gr}^{\mbox{\rm \scriptsize ad}}$
introduced (for a different 
purpose) in \cite{W1}. Let $V = \bigcap V_{\lambda}$ be the 
intersection of the $\lambda$-primary subspaces $V_{\lambda}$\,; 
let $k_{\lambda}$ be the codimension (in $\mathbb{C}[z]$) of 
$V_{\lambda}$\,, and set 
\smallskip
\begin{equation}
\label{adpoly}
m_V(z) = \prod_{\lambda}(z - \lambda)^{k_{\lambda}} \ .
\end{equation}
Let $W = m_V^{-1} V$\,, so that $W$ is a linear subspace of 
$\mathbb{C}(z)$\,. By definition, $\mbox{\rm Gr}^{\mbox{\rm \scriptsize ad}}$ 
consists of all subspaces $W \subset \mathbb{C}(z)$ that arise in 
this way. If we replace $V$ by $\rho V$  for some polynomial $\rho$\,, 
then $m_V$ gets replaced by $\rho m_V$\,; hence indeed the map 
$V \mapsto W$ induces a bijection from equivalence classes of primary 
decomposable subspaces to $\mbox{\rm Gr}^{\mbox{\rm \scriptsize ad}}$.
Note that the action of the group $\Gamma$ is compatible with the 
equivalence relation (\ref{equiv}), and so induces an action on 
$\mbox{\rm Gr}^{\mbox{\rm \scriptsize ad}}$. We write this action 
again as $W \mapsto \gamma W$\,, and think of it as scalar 
multiplication. If $U$ and $V$ are two points of 
$\mbox{\rm Gr}^{\mbox{\rm \scriptsize ad}}$, we define 
the set of differential operators $\mathcal{D}(U, V)$ as before 
(see (\ref{duv})).

Let us reformulate the Cannings-Holland results in terms of 
$\mbox{\rm Gr}^{\mbox{\rm \scriptsize ad}}$. If 
$W \in \mbox{\rm Gr}^{\mbox{\rm \scriptsize ad}}$, we set
$$
R_W = \mathcal{D}(\mathbb{C}[z], W) \equiv 
\{ D \in \mathbb{C}(z)[\partial] : D.\mathbb{C}[z]
\subset W \} \ .
$$
Clearly, $R_W$ is is a fractional right ideal of $A_1$ (we regard 
$\mathbb{C}(z)[\partial]$ as a subalgebra of the quotient field of 
$A_1 \equiv \mathbb{C}[z][\partial] $). 
Let $\alpha(W) \in \mathcal{R}$ be the 
isomorphism class of $R_W$\,; then Corollary \ref{chcor} 
implies immediately
\begin{theorem}
The map $\alpha : \mbox{\rm Gr}^{\mbox{\rm \scriptsize ad}}
\to \mathcal{R}$ is bijective; that is, every non-zero right ideal of $A_1$
is isomorphic to a unique ideal of the form $R_W$. We can
recover $W$ from $R_W$ by the formula
\begin{equation}
\label{getW}
W = e(R_W) \equiv \{ D.1 : D \in R_W \}\ .
\end{equation}
\end{theorem}

It is clear that the result of Corollary \ref{hardlem} 
is still valid for $U, V \in \mbox{\rm Gr}^{\mbox{\rm \scriptsize ad}}$. 
Applying it with $U = \mathbb{C}[z]$\,, $V = W$, we find that
$$
R_{\gamma W} = \gamma R_W \gamma^{-1} \ \, 
\text{for any} \ \, W \in \mbox{\rm Gr}^{\mbox{\rm \scriptsize ad}}\ ,
\ \gamma \in \Gamma\ .
$$
In other words, setting $\gamma(z) = e^{p(z)}$\,, we have 
\begin{proposition}
\label{phiprop}
Under the bijection $\alpha$, the action of 
the automorphism $\Phi_p$ on $\mathcal{R}$
corresponds to the map $W \mapsto e^{p(z)}W$ on
$\mbox{\rm Gr}^{\mbox{\rm \scriptsize ad}}$.
\end{proposition}

We can also assign to each point $W$ of
$\mbox{\rm Gr}^{\mbox{\rm \scriptsize ad}}$ the
fractional {\it left} ideal
$$
L_W = \mathcal{D}(W, \mathbb{C}[z]) \equiv
\{D \in \mathbb{C}(z)[\partial] : D.W \subset
\mathbb{C}[z] \} \ .
$$
This ideal is useful for studying the action of 
anti-automorphisms of $A_1$ on 
$\mbox{\rm Gr}^{\mbox{\rm \scriptsize ad}}$, because we have 
\begin{lemma}
\label{duallem}
The ideals $R_W$ and $L_W$ are dual to each other; that is,
in the notation of the previous section, we have
$ L_W = (R_W)' $\,.
\end{lemma}
\begin{proof}
It is obvious that $L_W \subset (R_W)'$\,. The converse
is a consequence of (\ref{getW}). Let
$w \in W$\,: then $w = D.1$ for some $D \in R_W$\,, so if
$q \in (R_W)'$\,, then $q.w = (qD).1$ is a polynomial, since
$qD \in A_1$\,. Hence $q \in L_W$\,.
\end{proof}
\section{The Baker function}

Associated to each point $ W $ of 
$ \mbox{\rm Gr}^{\mbox{\rm \scriptsize ad}} $ 
is the {\it Baker function} $ \psi_{W} $. This function plays an 
indispensable role, because both the bispectral involution
$ b $ on $ \mbox{\rm Gr}^{\mbox{\rm \scriptsize ad}} $
and the map $ \beta $ in (\ref{1.3}) are defined in term of
 $ \psi_{W} $. In this section we review briefly the definition
and the main properties of $ \psi_{W} $; the reader may 
consult \cite{SW}, \cite{W1} or \cite{W2} for more details.

Let $R_{-}$ be the space of rational functions that vanish
at infinity, so that we have
${\mathbb C}(z) = {\mathbb C}[z] \oplus R_{-}$\,.
If $W \in \mbox{\rm Gr}^{\mbox{\rm \scriptsize ad}}$\,, then 
the projection $\pi_W : W \to {\mathbb C}[z]$ defined by this
splitting is an operator of index zero (that is, its kernel and
cokernel have the same (finite) dimension); 
we say that $W$ belongs to the {\it big cell} of
$\mbox{\rm Gr}^{\mbox{\rm \scriptsize ad}}$ if $\pi_W$ is an
isomorphism. It is intuitively clear that ``most'' spaces $W$ 
belong to the big cell, in particular, that 
for any $W$, the space $\gamma W$ should belong to the
big cell for almost all $\gamma \in \Gamma$. We need the 
following more precise result, which is a consequence of 
Lemma 8.6 in \cite{SW}.
\begin{proposition}
\label{sw}
For any $W \in \mbox{\rm Gr}^{\mbox{\rm \scriptsize ad}}$, 
the space $e^{xz} W$ belongs to the big cell for all but a 
discrete\footnote{In fact it follows from Proposition \ref{to0} 
below that this set is finite}
set of values of $x \in \mathbb{C}$\,.
\end{proposition}

It follows that for any $W$ there is a unique function $\tilde \psi_W(\gamma, z)$ (defined for
almost all $\gamma \in \Gamma\,, \ z \in \mathbb{C}$) such that
\smallskip
\begin{verse}
(i) $\tilde \psi_W(\gamma, z)$ has the form $1 + O(z^{-1})$ for all
$\gamma \in \Gamma$\,;

(ii) $\tilde \psi_W(\gamma, z) \in \gamma^{-1} W$ for all
$\gamma \in \Gamma$\,.
\end{verse}
\smallskip
In (i) we have introduced two abbreviations of a kind that will 
be used without comment from now on. By ``$O(z^{-1})$'' we mean 
``some function of $z$ and possibly some other variables that 
vanishes as \,$z \to \infty$\,''; and ``for all $\gamma$'' means 
``for all $\gamma$ for which the statement makes sense'' 
(in the present case: for all $\gamma$ such that $\gamma^{-1} W$ 
belongs to the big cell).
The {\it Baker function} of $W$ is the function
$$
\psi_W(\gamma, z) = \gamma(z) \tilde \psi_W(\gamma, z) \ .
$$
We call $\tilde \psi_W$ the {\it reduced Baker function} of $W$.
The action of $\Gamma$ on $\mbox{\rm Gr}^{\mbox{\rm \scriptsize ad}}$ 
corresponds to translation in the first variable of the reduced 
Baker function; more precisely, if $\eta \in \Gamma$ and 
$W \in \mbox{\rm Gr}^{\mbox{\rm \scriptsize ad}}$, then 
\begin{equation}
\label{trans}
\tilde \psi_{\eta W}(\gamma ,z) = \tilde \psi_W(\eta^{-1} \gamma ,z)
\end{equation}
(for both sides have the form $1 + O(z^{-1})$ and belong to
$\eta \gamma^{-1} W$ for all $\gamma$).

Proposition \ref{sw} shows that we can restrict 
$\psi_W$ to the $1$-parameter subgroup
$\{ e^{xz} \}$ of $\Gamma$ : the resulting function of $x$ and $z$ 
is called the {\it stationary Baker function} of $W$\,;
we write it simply as $\psi_W(x,z)$ (rather than $\psi_W(e^{xz},z)$).
We refer to \cite{W1} for the proof of the 
following.
\begin{proposition}
\label{to0}
For any $W \in \mbox{\rm Gr}^{\mbox{\rm \scriptsize ad}}$, 
$\tilde \psi_W(x,z)$ is a separable rational function of $x$ 
and $z$ that tends to $1$ as either $x \to \infty$ or 
$z \to \infty$\,. 
\end{proposition}
Here ``separable'' means that $\tilde \psi_W(x,z)$ is a sum of 
functions of the form $f(x) g(z)$\,. 
We recall at this point that the {\it bispectral involution} $b$ on
$ \mbox{\rm Gr}^{\mbox{\rm \scriptsize ad}} $ is defined by the 
formula
\begin{equation}
\label{bisp}
\psi_{b(W)}(x,z) = \psi_{W}(z,x) \ \ .
\end{equation}
Of course, it is not {\it a priori} clear that the right hand side
of (\ref{bisp}) is the Baker function of any point of
$\mbox{\rm Gr}^{\mbox{\rm \scriptsize ad}}$; that this is the case
was the main result of \cite{W1}.

\begin{lemma}
\label{xlem}
Let $W \in \mbox{\rm Gr}^{\mbox{\rm \scriptsize ad}}$, and let
$\Theta(x)$ be any differential operator in $x$ with coefficients
that are (say) analytic on some open set in $\mathbb{C}$\,. Then
$$
(e^{-xz}\Theta(x)).\psi_W(x,z) \in e^{-xz} W \
\,\text{\rm for all} \,\ x\,.
$$
\end{lemma}
After multiplication by the polynomial $ m_{W} $ in
(\ref{adpoly}), this becomes an easy exercise in the space
$ {\mathcal H} $ of entire functions: we leave the details
to the reader.

The following consequence of 
Lemma \ref{xlem} is very useful.
\begin{corollary}
\label{useful}
Fix any $x$ such that the function $\psi_W(x,z)$ is regular 
at $x$\,. Then $e^{-xz} W$ is spanned by the functions
\begin{equation}
\label{i}
e^{-xz} \{ (\partial/\partial x)^i \psi_W(x,z) \}\ , \ \ i \geq 0 \ .
\end{equation}
\end{corollary}
\begin{proof}
By Lemma \ref{xlem}, the functions (\ref{i}) all belong to 
$e^{-xz} W$\,. They span $e^{-xz} W$ because the projection 
$e^{-xz} W \to \mathbb{C}[z]$ is an isomorphism, and the 
$i$th function (\ref{i}) has the form $z^i +$ (lower terms) 
for large $z$\,.
\end{proof}

\section{The algebra $\mathcal{W}$}

Proposition \ref{to0} shows that for large $x$ and $z$, the 
function $\psi_W$ has a series expansion of the form
$$
\psi_W(x,z) = e^{xz} \{ 1 + \sum_{i,j = 1}^{\infty} \alpha_{ij} 
x^{-j} z^{-i} \} \ .
$$
Formally, we can write this equation as $\psi_W = K_W.e^{xz}$\,, 
where $K_W$ denotes the formal integral operator
$$
K_W = 1 + \sum_{i,j = 1}^{\infty} \alpha_{ij} 
x^{-j} {\partial}_{x}^{-i} 
$$
(from now on we write ${\partial}_{x} \equiv \partial/\partial x$).
Obviously, the operator $ K_{W} $ contains exactly the same 
information as the function $ \psi_{W}(x,z)$; however, our proof of
Theorem \ref{th4} uses these operators in an
essential way, so we shall review the necessary formalism.

We denote by $\mathcal{W}$ the algebra of all formal operators 
of the form
\begin{equation}
\label{dsum}
L(x) = \sum_{i = - \infty}^{N} \sum_{j = - \infty}^{M} {\alpha}_{ij}
x^j \partial_{x}^{i}
\end{equation}
The (associative)
multiplication on $\mathcal{W}$ is uniquely determined by the 
constraint that it should coincide with the usual one on 
the subalgebra $\mathbb{C}((x^{-1}))[\partial_x]$ of 
differential operators. The operators $K_W$ above are 
invertible elements of $\mathcal{W}$\,. 
Now let $\mathcal{F}$ be the space of formal ``functions'' of the form
$$
f(x,z) = e^{xz} \left (
\sum_{i = - \infty}^{N} \sum_{j = - \infty}^{M} {\alpha}_{ij}
x^{j} z^{i} \right )\ .
$$
The {\it differential} operators in $\mathcal{W}$ 
act on $\mathcal{F}$ in an obvious way, and it is easy to check that
the action of $\partial_x$ is invertible, so that we can define the
action of $\mathcal{W}$ on  $\mathcal{F}$\,. In
this way $\mathcal{F}$ becomes a free
$\mathcal{W}$-module of rank 1
(as generator we can take $e^{xz}$).

On $\mathcal{W}$ we have the anti-automorphism $b$ defined by
$$
b(x) = \partial_x \ , \quad b(\partial_x) = x \ .
$$
On the other hand, because of the symmetry between $x$ and $z$ in
the space $\mathcal{F}$\,, the algebra
$\mathcal{W}(z)$ of formal 
operators in the variable $z$ also acts on
$\mathcal{F}$\,. If $L(x)$ is an operator of the form (\ref{dsum}),
we write $L(z)$ for the operator obtained from it by replacing $x$ by
$z$ and $\partial_x$ by $\partial \equiv \partial/\partial z$\,. 
We then have the following simple rule of calculation.
\begin{lemma}
\label{rule}
For any operator $L(x) \in \mathcal{W}$\,, we have
$$
L(x).e^{xz} = b(L)(z).e^{xz} \ .
$$
\end{lemma}

As a special case of Lemma \ref{rule}, we can reformulate
the definition (\ref{bisp}) of the bispectral involution
in terms of the operator $ K_{W} $. 

\begin{corollary}
\label{b}
The bispectral 
involution on $\mbox{\rm Gr}^{\mbox{\rm \scriptsize ad}}$ is 
characterized by the formula
$$
K_{b(W)} = b(K_W) \ .
$$
\end{corollary}
\begin{remark}
In the theory of integrable systems it is customary to work 
with rings of formal pseudo-differential operators in the 
style of Schur (see \cite{Sch}). Readers familiar with this
formalism should note that our algebra $ {\mathcal W} $ is
smaller than the ring 
$ {\mathbb C}((x^{-1}))((\partial_{x}^{-1})) $ of formal 
pseudo-differential operators with coefficients in
$ {\mathbb C}((x^{-1})) $; for example, the element
$ \sum_{i=1}^{\infty} x^{i} \partial_{x}^{-i} $ of this ring
does not belong to $ {\mathcal W} $. As this example shows, 
the involution $ b $ of $ {\mathcal W} $ cannot be extended
to the larger algebra: that is the main reason why we 
restrict our considerations to  $ {\mathcal W} $.  
\end{remark}

\section{Characterization of $\mathcal{D}(U,V)$}
In this section we work with the algebras $\mathcal{W}$ and 
$\mathcal{W}(z) $\,; we 
regard $\mathbb{C}(z)[\partial]$ as embedded in $\mathcal{W}(z)$ 
by identifying a rational function with its Laurent expansion 
near $ z = \infty$\,. 
\begin{lemma}
\label{+lem}
Let $U,\, V \in \mbox{\rm Gr}^{\mbox{\rm \scriptsize ad}}$\,,
$D \in \mathbb{C}(z)[\partial]$\,. Then there is a unique
differential operator $\Theta \in \mathbb{C}(x)[\partial_x]$
such that
\begin{equation}
\label{Oeq}
D(z).\psi_{U}(x,z) - \Theta(x).\psi_{V}(x,z) =
e^{xz} \{O(z^{-1}) \} \ .
\end{equation}
Explicitly, we have $\Theta = [K_U b(D)(x) K_{V}^{-1}]_{+}$\,, 
where the subscript {\rm +} means that we delete the terms involving 
negative powers of $\partial_x$\,. 
\end{lemma}
\begin{proof}
Note first that
$$
D(z).\psi_{U}(x,z) = D(z) K_U(x).e^{xz} = K_U(x) D(z).e^{xz} =
K_U(x) b(D)(x).e^{xz}
$$
(we used the rule in Lemma \ref{rule} at the last step). So the
formula (\ref{Oeq}) is equivalent to
$$
[K_U b(D)(x) - \Theta K_V].e^{xz} = e^{xz} \{ O(z^{-1}) \} \ .
$$
That is true if and only if the operator acting
on $e^{xz}$ on the left has negative order; equivalently,
if the operator $K_U b(D)(x) K_{V}^{-1} - \Theta$ has negative
order. If $\Theta$ is differential, this means exactly that
$\Theta = [K_U b(D)(x) K_{V}^{-1}]_{+}$\,. Since $K_U$ and $K_V$
(and hence also $K_{V}^{-1}$) have rational coefficients, and
$b(D)$ has polynomial coefficients, it is clear that
$\Theta \in \mathbb{C}(x)[\partial_x]$\,.
\end{proof}
\begin{proposition}
\label{keylem}
Let $U,\, V \in \mbox{\rm Gr}^{\mbox{\rm \scriptsize ad}}$,
$D \in \mathbb{C}(z)[\partial]$\,. Then the following are equivalent: \\
{\rm (i)} $D \in \mathcal{D}(U, V)$\ ; \\
{\rm (ii)} There is a differential operator $\Theta(x)$ such that
$D(z).\psi_U(x,z) = \Theta(x).\psi_V(x,z)$\,; \\
{\rm (iii)} The operator $K_U b(D)(x) K_V^{-1}$ is differential.
\end{proposition}
\begin{proof}
The equivalence of (ii) and (iii) follows at once from
Lemma \ref{+lem}. To prove the equivalence of (i) and
(ii), we have only to recycle some standard arguments
(see, for example \cite{Kri}, \cite{SW}).
First, suppose that $D.U \subset V$\,, and let $\Theta$ be
the unique differential operator such that (\ref{Oeq})
holds, or, equivalently, such that
\begin{equation}
\label{EOeq}
(e^{-xz} D(z) e^{xz}).\tilde{\psi}_U(x,z) -
(e^{-xz} \Theta(x)).{\psi}_V(x,z) = O(z^{-1}) \ .
\end{equation}
By Proposition \ref{phiprop}, the first term here belongs to $e^{-xz} V$ for
all $x$\,; by Lemma \ref{xlem}, so does the second term. Since
generically $e^{-xz} V$ contains no function of the form
$O(z^{-1})$\,, it follows that the expression on the left
of (\ref{EOeq}) is zero, that is, $D.\psi_U = \Theta.\psi_V$\,.
Conversely, suppose $D.\psi_U = \Theta.\psi_V$\,.
Differentiating with respect to $x$\,, we get
\begin{equation}
\label{mess}
\{ e^{-xz} D(z) e^{xz} \}.\{ e^{-xz} (\partial_{x}^{i}.\psi_U(x,z)) \} =
e^{-xz} (\partial_{x}^{i} \Theta(x)). \psi_V(x,z) \ \,
\end{equation}
for all $i \geq 0$ and for all $x$\,. 
Fix any $x$ such that the functions $\psi_U(x,z)$ and 
$\psi_V(x,z)$ are regular at $x$\,. By Lemma \ref{xlem}, 
the right hand side of (\ref{mess}) belongs to $e^{-xz} V$\,. 
By Corollary \ref{useful}, the functions 
$e^{-xz} (\partial_{x}^{i}.\psi_U(x,z))$  
on the left hand side 
of (\ref{mess}) span $e^{-xz} U$\,; hence the operator 
$e^{-xz} D(z) e^{xz}$ maps $e^{-xz} U$ into $e^{-xz} V$\,. 
By Corollary \ref{hardlem} we then have $D.U \subset V$\,.
\end{proof}

We shall use Proposition \ref{keylem} in the following form.
\begin{corollary}
\label{main}
Let $U,\, V \in \mbox{\rm Gr}^{\mbox{\rm \scriptsize ad}}$,
$D \in \mathcal{W}(z)$\,. Then the following are equivalent: \\
{\rm(i)} $D \in \mathcal{D}(U, V)$\ ; \\
{\rm (ii)} the operators $D$ and $K_U(z) b(D) K_V^{-1}(z)$ 
are both differential.
\end{corollary}
\begin{proof}
The only thing left to see is that condition (ii) ensures that 
the differential operator $D \in \mathcal{W}(z)$ has rational coefficients.
The argument is like that in the proof of Lemma \ref{+lem} above:
if we set $\Theta = K_U(z) b(D) K_V^{-1}(z)$ then we have
$$
D = b(K_U)^{-1}(z) b(\Theta) b(K_V)(z) \ .
$$
Because the Baker functions $\tilde \psi_W$ are rational in both 
variables $x$ and $z$, the operators $b(K_W)$, hence also $b(K_W)^{-1}$, 
have rational coefficients; and because $\Theta$ is differential, 
$b(\Theta)$ has polynomial coefficients. It follows that $D$ has 
rational coefficients.
\end{proof}
\section{The action of $\widehat G$ on the Grassmannian}

As we observed in the Introduction, the action of the automorphisms 
$\Phi_p$ and of the anti-automorphisms $b$ and $c$ in (\ref{ab}) correspond 
to known symmetries of $\mbox{\rm Gr}^{\mbox{\rm \scriptsize ad}}$. 
We have already seen this in the case of the $\Phi_p$ 
(see Proposition \ref{phiprop}). 

Next, we prove the crucial Theorem \ref{th4} from the Introduction.
Using Lemma \ref{duallem} and the
definition (\ref{act2}), we see that Theorem \ref{th4}
amounts to the following assertion.
\begin{proposition}
\label{bprop}
For any $W \in \mbox{\rm Gr}^{\mbox{\rm \scriptsize ad}}$,
the right ideals $R_{b(W)}$ and $b(L_W)$ are isomorphic.
\end{proposition}

The (fractional) ideals $ R_{b(W)} $ and $ b(L_{W}) $ are contained
(respectively) in the subalgebras $ {\mathbb C}(z)[\partial] $
and $ {\mathbb C}(\partial)[z] $ of the quotient field $ Q $
of $ A_1 $. Since these algebras are both canonically embedded
in $ {\mathcal W}(z) $ (by replacing rational functions by their Laurent 
expansions near infinity), we can prove Proposition~\ref{bprop}
by calculating in the algebra $ {\mathcal W}(z) $ (rather than in $ Q $).
We shall prove the following more precise result.
\begin{proposition}
\label{Kprop}
For any $W \in \mbox{\rm Gr}^{\mbox{\rm \scriptsize ad}}$, we have
$$
R_{b(W)} = K_{W}(z) b(L_W) \ .
$$
\end{proposition}
\begin{proof}
In the proof all operators are supposed to be in the variable 
$z$\,: we omit $z$ from the notation. 
Applying Corollary \ref{main} with $U = \mathbb{C}[z]$\,,
$V = b(W)$\,, we get that
\begin{equation}
\label{Req}
D \in R_{b(W)} \Longleftrightarrow D \ \,\mbox{and} \ \,
b(D)(K_{b(W)})^{-1}
\ \,\text{\rm are differential} \ .
\end{equation}
Also, $D \in K_{W} b(L_W)$ if and only if
$b(D)b(K_W)^{-1} \in L_W$\,, so
applying Corollary \ref{main} again 
with $U = W$\,, $V = \mathbb{C}[z]$\,,
we get
\begin{equation}
\label{Leq}
D \in K_Wb(L_W) \Longleftrightarrow b(D)b(K_W)^{-1}
\ \,\mbox{and} \ \,K_W[(K_W)^{-1}D] \,\ \text{\rm are differential} \ .
\end{equation}                       
But in view of Corollary~\ref{b}, the conditions~(\ref{Req}) and   
(\ref{Leq}) coincide; hence the Proposition.
\end{proof}

To end this section, we identify the action of the formal adjoint $c$ 
on $\mbox{\rm Gr}^{\mbox{\rm \scriptsize ad}}$. Because the 
group $\widehat G$ is generated by the automorphisms 
$\Phi_p$ and the anti-automorphism $b$\,, this result is 
actually surplus to our needs; however, it seems of some 
interest in its own right.
\begin{proposition}
The action of $c$ on $\mbox{\rm Gr}^{\mbox{\rm \scriptsize ad}}$
is given by $c(W) = W^*$, where (as in \cite{W2}) $W^*$
denotes the annihilator of $W$ with respect to the symmetric
bilinear form
$$
\langle f,g \rangle = \mbox{\rm res}_{\infty}\, f(z)g(z)dz
$$
on $\mathbb{C}(z)$\,.
\end{proposition}
\begin{proof}
According to the definitions, we have to show that
$R_{W^*}$ is isomorphic to $c(L_W)$\,. In fact these two 
right ideals coincide, for we have (denoting the formal adjoint 
of an operator $D$ by $D^*$)
\begin{eqnarray}
D \in R_{W^*} & \Longleftrightarrow &
Df \in W^* \ \ \mbox{\rm for\ all} \  \, f \in \mathbb{C}[z] \nonumber \\
& \Longleftrightarrow & \langle D.f, g \rangle  = 0 \ \
\mbox{\rm for\ all} \  \,f \in \mathbb{C}[z]\ ,\ g \in W \nonumber \\
& \Longleftrightarrow & \langle f, D^*.g \rangle  = 0 \ \
\mbox{\rm for\ all} \  \, f \in \mathbb{C}[z]\ ,\ g \in W \nonumber \\
& \Longleftrightarrow & D^*g \in \mathbb{C}[z] \
\ \mbox{\rm for\ all} \  \, g \in W  \nonumber \\
& \Longleftrightarrow & D^* \in L_W  \nonumber \\
& \Longleftrightarrow & D \in ({L_W})^* \equiv c(L_W) \ . \nonumber
\end{eqnarray}
\end{proof}
\section{Proof of Theorems \ref{th2} and \ref{th3}}

It is now easy to see how
the action of $\widehat G$ on $\mbox{\rm Gr}^{\mbox{\rm \scriptsize ad}}$
transfers to $\mathcal{C}$\,. First, the formula (\ref{phiEq}) for
the action of $\Phi_p$ is just the basic fact (see \cite{W2}) that the
Calogero-Moser flows and the KP flows correspond under $\beta$\,.
Indeed, according to the 
definition\footnote{After the change of notation 
$X \mapsto -X^t, \ Z \mapsto -Y^t$\,.}
in \cite{W2}, $\beta$ maps (the conjugacy class of) a pair of 
matrices $(X,Y)$ to the point $W$  with reduced Baker function
\begin{equation}
\label{bak}
\tilde \psi_W(\gamma ,z) = \det \left\{ \mbox{\rm I} -
(p'(Y) - X)^{-1}(z \mbox{\rm I} - Y)^{-1} \right \} \ ,
\end{equation}
where $\gamma(z) = e^{p(z)} \in \Gamma$\,.
If now $q(z)$ is any polynomial with zero constant term,
and we set $\eta (z) = e^{q(z)}$, then the formula (\ref{trans}) 
shows that $\tilde \psi_{\eta W}$
is obtained by replacing $X$ by $X + q'(Y)$ in the formula
(\ref{bak}); that is, we have
$$
\beta(X + q'(Y),\, Y) = e^{q(z)} W \ ,
$$
which is what the formula (\ref{phiEq}) states (after 
restoring the notation $ p $ for $ q $).

On the other hand, the bispectral involution $b$ corresponds
under $\beta$ to the map $(X,Y) \mapsto (Y^t, X^t)$ on pairs
of matrices.
So the formula (\ref{psiEq}) follows from (\ref{phiEq})
and the fact that $\Psi_q = b \,\Phi_{-q} b$\,. That completes 
the proof of Theorem \ref{th2}.

For the proof of Theorem \ref{th3}, 
we choose as base-point in
${\mathcal C}_{n}$ the pair $(X_0 , Y_0)$ given by
$$
X_0 = -\sum_{r=1}^{n-1} r E_{r+1,\,r}\quad , \qquad
Y_0 = \sum_{r=1}^{n-1} E_{r,\,r+1} \, 
$$
(as usual, $E_{i,j}$ denotes the matrix with $(i,j)$-entry $1$ 
and zeros elsewhere). We shall prove something a little more 
precise than Theorem \ref{th3}: namely, 
we shall show that any given point $(X,Y) \in {\mathcal C}_{n}$
can be moved to the base-point by applying at most three of
the automorphisms $\Phi_p$ and $\Psi_q$ (in fact two are
enough unless $X$ and $Y$ are both non-diagonalizable).
\begin{lemma}
\label{plem}
Let ${\boldsymbol \lambda} = \{\lambda_1, \ldots , \lambda_n\}$ be any
collection of \,$n$ complex numbers. Then there is a unique
polynomial $p(z)$ of degree $n$ and with zero constant term
such that the set of eigenvalues of the matrix
$X_0 + p'(Y_0)$ is ${\boldsymbol \lambda}$\,.
\end{lemma}
\begin{proof}
According to \cite{Mac} (Ch.\ 1, Sect.\ 2, Example 8),
$\det \lbrace X_0 + p'(Y_0) \rbrace$ is
equal to $n!$ times the coefficient of $z^n$ in
the series $\exp \lbrace p(z) \rbrace$\,. The lemma follows
easily from this.
\end{proof}
\begin{lemma}
\label{dlem}
Let $(X,Y) \in {\mathcal C}_{n}$\,, and suppose that $X$ is
diagonalizable. Then there are unique polynomials $p$ and $q$
of degree $n$ and with zero constant term
such that
$$
\Psi_q \Phi_p (X_0, Y_0) = (X, Y) \, .
$$
\end{lemma}
\begin{proof}
We may assume that $X = \mbox{\rm diag}(x_1 , \ldots , x_n)$
is diagonal and that $Y$ is a Calogero-Moser matrix,
that is, the off-diagonal entries of $Y$ are given by
$$
Y_{ij} = (x_i - x_j)^{-1} \ \, \text{for} \ \, i \neq j 
$$
(cf.\ \cite{W2}, (1.14): the $x_i$ are necessarily 
distinct). By Lemma \ref{plem}, there
is a unique $p$ such that $X_1 = X_0 + p'(Y_0)$
has the same eigenvalues as $X$; then $(X_1, Y_0)$ is conjugate
to a pair $(X, Y_1)$ with $Y_1$ another Calogero-Moser matrix.
Thus \,$Y$ and\, $Y_1$ differ only in their diagonal entries. 
From the non-vanishing of the Vandermonde determinant
$\det(x_i^j)$\,, we find
that there is a unique polynomial $q$ as in the Lemma
such that $Y = Y_1 - q'(X)$. The
Lemma follows.
\end{proof}
\begin{lemma}
\label{taka}
Let $(X,Y) \in {\mathcal C}_{n}$\,. Then there is a polynomial
$r$ such that $X + r'(Y)$ is diagonalizable.
\end{lemma}
\begin{proof}
This follows from Shiota's lemma (see \cite{W2}, Lemma 5.6).
\end{proof}
{\it Proof of Theorem \ref{th3}}. The Theorem now follows at
once from Lemmas \ref{dlem} and \ref{taka}.
\section{$G$ as an algebraic group}
Nakajima (see \cite{N}, \cite{W2}) 
has found that the Hilbert scheme $\mbox{\rm Hilb}_n(\mathbb{A}^2)$ 
is a hyperk\"ahler manifold, and that after changing the complex 
structure to a different one in the hyperk\"ahler family, we 
obtain the space $\mathcal{C}_n$\,. It is difficult to resist 
the feeling that this deformation of complex structure from 
$\mbox{\rm Hilb}_n(\mathbb{A}^2)$ to $\mathcal{C}_n$ is related 
to the deformation of algebras from $A_0$ to $A_1$ (via the 
algebras $A_{\lambda}$ with commutation relation 
$[x,y] = \lambda$). It would be interesting to understand this analogy 
more precisely; here we just want to point out that considering 
the action of the group $G$ on $\mathcal{C}_n$  
might lead to a new perspective on Nakajima's result. 
 
So far we have considered $G$ simply as an abstract (discrete)
group. However, the fact that it acts transitively on the affine 
algebraic varieties $\mathcal{C}_n$ suggests that we should try to 
view $G$ as an {\it algebraic} group in such a way that its 
action on $\mathcal{C}_n$ is algebraic. Further,
from the formulae (\ref{phiEq}), 
(\ref{psiEq}), it is clear that the action of $G$ preserves the 
natural holomorphic symplectic structure on the spaces 
$\mathcal{C}_n$ (see \cite{W2}, p.\ 9): thus 
$\mathcal{C}_n$ 
should be a coadjoint orbit of $G$ (or possibly of some central 
extension of $G$). 
The fact that these orbits are hyperk\"ahler could then perhaps
be compared to Kronheimer's result that the coadjoint orbits
of finite dimensional complex semisimple groups have hyperk\"ahler
structures (see \cite{Kr}, \cite{Ko}, \cite{B}). As N. Hitchin
has pointed out to us, the phemonenon that the transitive action
of $G$ on $\mathcal{C}_n$ becomes intransitive on
$\mbox{\rm Hilb}_n(\mathbb{A}^2)$ then falls into perspective, 
since it is just what happens in Kronheimer's case.  

Let us recall first that $G$ can be identified with the group
$G_0$ of unimodular\footnote{that is, with Jacobian determinant 1} 
automorphisms of  the commutative algebra 
$A_0 \equiv \mathbb{C}[x,y]$\,. Indeed, it is well known (see \cite{J}) that 
$G_0$ is generated by the automorphisms $\Phi_p$ and $\Psi_q$ 
defined by the formulae (\ref{phieq}) and (\ref{psieq}); 
furthermore (see \cite{M2}) the relations between these are 
the same in $G_0$ as in $G$\,, so that as abstract groups we 
may identify $G_0$ and $G$\,. 
The group $G_0$ is a prototype example in Shafarevich's theory 
of infinite dimensional algebraic groups (see \cite{Sh}). The 
key idea is that we give $G_0$ a structure of infinite dimensional algebraic 
variety by regarding it as the union of the finite dimensional 
subvarieties $G^{(d)}$ of automorphisms of degree at most $d$\,: here 
the ``degree'' of an automorphism $\sigma$ is defined to 
be the largest of the degrees of the polynomials 
$$
\sigma(x)\,,\ \sigma(y)\,;\ \sigma^{-1}(x)\,,\ \sigma^{-1}(y)\ .
$$
This idea works equally well for the group $G$ of automorphisms 
of $A_1$\,: let us denote the resulting algebraic group by 
$G_1$\,.  
We claim that the algebraic groups $G_1$ and  
$G_0$ are not isomorphic. Indeed, Shafarevich shows in \cite{Sh} 
that the Lie algebra of $G_0$ (considered as an algebra of 
derivations of $A_0$) is the algebra
of all (polynomial) vector fields with zero divergence on
$\mathbb{A}^2$\,; that is, it is isomorphic to $A_0/{\mathbb{C}}$\,,
where $A_0$ is made into a Lie algebra via the canonical
Poisson bracket. A similar calculation for $G_1$ shows that its
Lie algebra is the full algebra of derivations of $A_1$\,, so
it is isomorphic to $A_1/{\mathbb{C}}$\,. Since their Lie algebras
are not isomorphic, the algebraic groups $G_0$ and $G_1$ are not
isomorphic. Thus as we
deform $A_0$ to $A_1$\,, the complex structure on the space of
ideals deforms from $\mbox{\rm Hilb}_n(\mathbb{A}^2)$
to $\mathcal{C}_n$\,, and the automorphism group deforms from
$G_0$ to $G_1$\,. Certainly, $G_0$ acts algebraically on 
$\mbox{\rm Hilb}_n(\mathbb{A}^2)$, so it 
is very natural to expect at this point that 
$G_1$ should act algebraically on $\mathcal{C}_n$\,; however, that is not 
the case. In fact the group $G_1$ has no non-trivial finite dimensional
homogeneous spaces; for if it did, we should obtain a non-trivial
homomorphism from $A_1/{\mathbb{C}}$ to the Lie algebra of vector
fields on a finite dimensional variety, and it is known that
this cannot happen (see \cite{HM}, \cite{St1}). In a similar way, 
we have convinced ourselves that for $n > 1$ the action of $G_0$ on 
$\mathcal{C}_n$ is not algebraic (despite an implied belief 
to the contrary in \cite{A}).

These are discouraging facts; 
however, there is yet a third structure of algebraic group on $G$ 
that does seem to answer our purpose: namely, we can identify $G$ with 
a group of automorphisms of the free associative algebra 
$A = \mathbb{C} \langle x, y \rangle$\,. The abelianization map 
$A \to A_0$ induces a natural homomorphism from $\mbox{\rm Aut}(A)$ 
to $\mbox{\rm Aut}(A_0)$, and it is known that this homomorphism 
is bijective (see \cite{Cz}, \cite{M1}). Denoting now by 
$\mathcal{G}$ the subgroup of ``unimodular'' 
(that is, preserving the commutator $xy - yx$) automorphisms 
of $A$\,, we obtain bijective maps of algebraic groups 
$\mathcal{G} \to G_0$ and 
$\mathcal{G} \to G_1$\,. We claim that neither 
of these maps is an isomorphism; in fact the induced maps of 
Lie algebras are not injective. Thus $G_0$ and $G_1$ should be 
thought of as distinct quotient groups of $\mathcal{G}$\,.

It seems to us that $\mathcal{G}$ does indeed act 
algebraically on the varieties $\mathcal{C}_n$\,. We think that these 
varieties deserve further study from this point of view\footnote{After 
reading an earlier version of the present paper, V. Ginzburg 
(see \cite{G}) has shown that $\mathcal{C}_n$ can indeed be 
identified with a coadjoint orbit of a central extension of $\mathcal{G}$}.

\bibliographystyle{amsalpha}

\end{document}